\documentclass{article}
\usepackage[utf8]{inputenc}
\usepackage[pass]{geometry}

\usepackage[dvipsnames]{xcolor}

\usepackage[utf8]{inputenc}
\usepackage[T1]{fontenc}
\usepackage{lmodern}

\usepackage{listings}

\usepackage{amsmath}
\usepackage{mathtools}
\usepackage{amssymb}
\usepackage{amsthm}

\usepackage{array}
\usepackage{comment}

\usepackage[plainpages=false,pdfpagelabels,colorlinks=true,citecolor=blue,hypertexnames=false]{hyperref}
                       
\newmuskip\pFqmuskip

\usepackage{todonotes}

\makeatletter
\def\@fnsymbol#1{\ensuremath{\ifcase#1\or \star\or \ddagger\or
   \mathsection\or \mathparagraph\or \|\or **\or \dagger\dagger
   \or \ddagger\ddagger \else\@ctrerr\fi}}
\makeatother

\newcommand*\pFq[6][8]{%
  \begingroup 
  \pFqmuskip=#1mu\relax
  \mathcode`\,=\string"8000
  \begingroup\lccode`\~=`\,
  \lowercase{\endgroup\let~}\pFqcomma
  {}_{#2}F_{#3}{\left[\genfrac..{0pt}{}{#4}{#5};#6\right]}%
  \endgroup
}     
\newcommand{\pFqcomma}{\mskip\pFqmuskip}

\newtheorem{thm}{Theorem}
\newtheorem{prop}{Proposition}

\newtheorem{lem}{Lemma}

\newtheoremstyle{note}
{3pt}
{3pt}
{}
{}
{\bfseries}
{.}
{.5em}
{}
\theoremstyle{note}
\newtheorem{rem}{Remark}

\newcommand{\Iso}{\textrm{\sf Iso}}

\title{A hypergeometric proof that $\Iso$ is bijective}
\author{Alin Bostan\footnote{Inria, Univ. Paris-Saclay, France, \url{alin.bostan@inria.fr}.} 
\ and Sergey Yurkevich\footnote{U. Wien, Austria and Inria, Univ. Paris-Saclay, France, \url{sergey.yurkevich@univie.ac.at}.}}
\date{\today}
 
\begin{document}

\maketitle    

\begin{abstract}
We provide a short and elementary proof of 
the main technical result of the recent article ``\emph{On the uniqueness of Clifford torus with prescribed isoperimetric ratio}''~\cite{YuCh20} by Thomas Yu and Jingmin Chen. The key of the new proof is an explicit expression of the central function ($\Iso$, to be proved 
bijective)
as a quotient of Gaussian hypergeometric functions.       
\end{abstract}	                                         

In their recent paper~\cite{YuCh20}, 
Thomas Yu and Jingmin Chen needed to prove,
as a crucial intermediate result, 
that a certain real-valued function $\Iso$
(related to isoperimetric ratios of Clifford tori) is monotonic increasing.   
They reduced the proof of this fact to 
the positivity of a sequence of rational numbers $(d_n)_{n \geq 0}$, defined explicitly in terms of nested binomial sums. 
This positivity was subsequently proved
by Stephen Melczer and Marc Mezzarobba~\cite{MeMe20},
who used a computer-assisted approach
relying on analytic combinatorics and rigorous numerics, combined with the fact 
(proved in~\cite{YuCh20}) that the sequence 
    $(d_n)_{n\geq 0}$ satisfies an explicit linear recurrence of order seven with polynomial coefficients in $n$.

In this note, we provide an alternative, short and conceptual, proof 
of the monotonicity of the function \Iso.
Our approach is different in spirit from the ones in~\cite{YuCh20} and~\cite{MeMe20}. 
Our 
main result (Theorem~\ref{thm:2F1-Iso} below) is that the function 
$\Iso(z)$ can be expressed in terms of Gaussian hypergeometric functions $_2F_1$ defined by       
\begin{equation}\label{eq:2F1}
\pFq{2}{1}{a,b}{c}{z} = \sum_{n=0}^\infty \frac{(a)_n(b)_n}{(c)_n} \, \frac
{z^n} {n!},
\end{equation}
where $(a)_n$ denotes the 
rising factorial $(a)_n=a(a+1)\cdots(a+n-1)$ for $n\in\mathbb{N}$.     

\smallskip In the notation of Yu and Chen, the function 
\[\Iso:[0,\sqrt{2}-1) \rightarrow [3/2 \cdot (2\pi^2)^{-1/4},1)\]
is given as
\begin{equation} \label{def:iso}
\Iso(z) = 6 \sqrt{\pi} \cdot \frac{V(z)}{A^{3/2}(z)},
\end{equation}
where $A(z) = \sum_{n \geq 0} a_n z^{2n}$ and   
$V(z) = \sum_{n \geq 0} v_n z^{2n}$ 
are complex analytic functions in the disk $\{ z : |z| < \sqrt{2}-1 \}$, 
given by the power series expansions 
\[
A(z) = \sqrt{2} \pi^2 \cdot \left(4 + 52 z^2 + 477 z^4 + 3809 z^6 + \frac{451625}{16} z^8 + \cdots \right),
\]
\[
V(z) = \sqrt{2} \pi^2 \cdot \left(2 + 48 z^2 + \frac{1269}{2}  z^4 + 6600 z^6 + \frac{1928025}{32} z^8 + \cdots \right).
\]   
The precise definitions of $A$ and $V$ are given in Section 4.3 of~\cite{YuCh20}, notably in equations (4.2)--(4.3). Since the sequences $(a_n)_{\geq 0}$ and $(v_n)_{\geq 0}$ are expressed in terms of nested binomial sums, $A(z)$ and $V(z)$ satisfy linear differential equations with polynomial coefficients in $z$, that can be found and proved automatically using \emph{creative telescoping}~\cite{Chyzak00}.
Yu and Chen, resp. Melczer and Mezzarobba, use this methodology to 
find 
a linear recurrence satisfied by the coefficients $(d_n)_{n\geq 0}$ of 
\[
F(z):= \frac{1}{4\pi^4} \cdot \left( \frac{2V'(\sqrt{z})A(\sqrt{z}) - 3 V(\sqrt{z}) A'(\sqrt{z})}{\sqrt{z}}\right)
= 72  + 1932 z + 31248 z^3 + \cdots,
\]
respectively a linear differential equation satisfied by the function $F(z)$.

\smallskip \noindent Similarly, one can compute linear differential equations satisfied individually by 
\[
\bar{A}(z):=\frac{1}{\sqrt{2} \pi^2} \cdot A(\sqrt{z}) = 
4+52\,z+477\,{z}^{2}+3809\,{z}^{3}+{\frac{451625}{16}}{z}^{4}+{\frac{
3195333}{16}}{z}^{5}+ \cdots
\]
and by
\[
\bar{V}(z):=\frac{1}{\sqrt{2} \pi^2} \cdot V(\sqrt{z}) = 
2+48\,z+{\frac{1269}{2}}{z}^{2}+6600\,{z}^{3}+{\frac{1928025}{32}}{z}
^{4}+{\frac{2026101}{4}}{z}^{5}+ \cdots.           
\]
Concretely, $\bar{A}(z)$ and $\bar{V}(z)$ %
satisfy second-order linear differential equations:
\begin{align*}
z ( z-1 ) ( {z}^{2}-6z+1 )  ( z+1) ^{2}
\bar{A}''(z)
+ 
& ( z+1 )  ( 5{z}^{4}-8{z}^{3}-32{z}^{2
}+28z-1 ) \bar{A}'(z)  
\\
&+\left( 4{z}^{4}+11{z}^{3}-{z}^{2}-43z+13 \right) \bar{A}(z)
= 0 
\end{align*}
and respectively
\begin{align*}
z \left( z-1 \right)  \left( z+1
 \right)  \left( {z}^{2}-6\,z+1 \right) ^{2}\bar{V}''(z)         
 & \\ 
+ \left( {z}^{2}-6\,z+1 \right)  \left( 7
\,{z}^{4}-22\,{z}^{3}-18\,{z}^{2}+26\,z-1 \right) \bar{V}'(z)
& \\  +   
3\, \left( 3\,{z}^{5}-24\,{z}^{4}-2\,{z}^{3}+56\,{z}^{2}-25\,z+8
 \right) \bar{V}(z)
& = 0.  
\end{align*}  

\noindent From these equations, we deduce the following closed-form expressions:
\begin{thm}  \label{thm:2F1-A-V}
The following equalities hold for all $z\in\mathbb{R}$ with $0 \leq z \leq \sqrt{2}-1$:
\[
\bar{A}(z) =   
\frac{4 \left(1-z^{2}\right)}{\left(z^{2}-6 z +1\right)^{2}}  
\, \cdot \, \pFq{2}{1}{-\frac{1}{2},-\frac{1}{2}}{1}{\frac{4 z}{\left(1-z \right)^{2}}}
\]
and
\[
\bar{V}(z) =
\frac{2 \left(1-z\right)^{3}}{\left(z^{2}-6 z +1\right)^{3}}  
\, \cdot \, \pFq{2}{1}{-\frac{3}{2},-\frac{3}{2}}{1}{\frac{4 z}{\left(1-z \right)^{2}}}       .
\]
\end{thm}

\begin{proof}  
It is enough to check that the right-hand side expressions  satisfy the same linear differential equations 
as $\bar{A}$ and $\bar{V}$, with the same initial conditions.
\end{proof}   
        
As a direct consequence of Theorem~\ref{thm:2F1-A-V} and of definition~\eqref{def:iso} we get:
\begin{thm}           \label{thm:2F1-Iso}
The function \emph{$\Iso$} admits the following closed-form expression:
\[
\emph{\Iso}^2(z) = 
\frac{9  \sqrt{2}}{8 \pi} \cdot
\frac{\pFq{2}{1}{-\frac{3}{2},-\frac{3}{2}}{1}{\frac{4 z^2}{\left(1-z^2 \right)^{2}}}^2}{\pFq{2}{1}{-\frac{1}{2},-\frac{1}{2}}{1}{\frac{4 z^2}{\left(1-z^2 \right)^{2}}}^3}    
\cdot \left(\frac{1-z^2}{1+z^2}\right)^3.
\]
\end{thm} 

Using the expression in Theorem~\ref{thm:2F1-Iso}, we can now prove the main result of~\cite{YuCh20}. 

\begin{thm}           \label{thm:mono-Iso}
\emph{$\Iso$} is a monotonic increasing function and $\lim_{z\rightarrow \sqrt{2}-1} \emph{\Iso}(z)=1$. In particular,  \emph{$\Iso$} is a bijection.
\end{thm} 
     
\begin{proof}
The value of ${\Iso}^2(z)$ at $z=\sqrt{2}-1$ is equal to 
\[
{\Iso}^2(\sqrt{2}-1) = 
\frac{9  \sqrt{2}}{8 \pi} \cdot
\frac{\pFq{2}{1}{-\frac{3}{2},-\frac{3}{2}}{1}{1}^2}{\pFq{2}{1}{-\frac{1}{2},-\frac{1}{2}}{1}{1}^3}    
\cdot \frac{\sqrt{2}}{4}.
\]    
From Gauss's summation theorem~\cite[Th.~2.2.2]{AnAsRo99} it follows that $\pFq{2}{1}{-\frac{3}{2},-\frac{3}{2}}{1}{1} = 32/(3 \pi)$ and
$\pFq{2}{1}{-\frac{1}{2},-\frac{1}{2}}{1}{1} = 4/\pi$; therefore,
\[
{\Iso}^2(\sqrt{2}-1) = 
\frac{9  \sqrt{2}}{8 \pi} \cdot
\frac{(32/(3 \pi))^2}{(4/\pi)^3}    
\cdot \frac{\sqrt{2}}{4} = 1.
\]          
It remains to prove that   {$\Iso$} is monotonic increasing.
It is enough to show that 
\[
z \; \mapsto \; \frac{\pFq{2}{1}{-\frac{3}{2},-\frac{3}{2}}{1}{\frac{4 z}{\left(1-z \right)^{2}}}^2}{\pFq{2}{1}{-\frac{1}{2},-\frac{1}{2}}{1}{\frac{4 z}{\left(1-z \right)^{2}}}^3}    
\cdot \left(\frac{1-z}{1+z}\right)^3
\]                                  
is increasing on $[0, 3-2\sqrt{2})$. Equivalently, via the change of variables
$ x = \frac{4 z}{\left(1-z \right)^{2}}$, it is enough to prove that the function 
\[
h \colon x \; \mapsto \; \frac{\pFq{2}{1}{-\frac{3}{2},-\frac{3}{2}}{1}{x}^2}{\pFq{2}{1}{-\frac{1}{2},-\frac{1}{2}}{1}{x}^3}    
\cdot  \left( x+1 \right) ^{-{\frac{3}{2}}}
\]                                  
is increasing on $[0, 1)$.   
Clearly, $h$ can be written as $h = f^3 \cdot g^2$, where
\[
f(x) 
=
\frac{\sqrt{x+1}}{\pFq{2}{1}{-\frac{1}{2}, -\frac{1}{2}}{1}{x}}    
\quad \text{and} \quad 
g(x) = \frac{\pFq{2}{1}{-\frac{3}{2}, -\frac{3}{2}}{1}{x}}{(x+1)^{\frac32}}     .
\]                 
Hence, it is enough to prove that both $f$ and $g$ are increasing on $[0, 1)$.  
We will actually prove a more general fact in Proposition~\ref{prop:1}, which may be of independent interest.   
Using that $w_{1/2} = 1/f$ and $w_{3/2} = g$, we deduce from Proposition~\ref{prop:1} that both $f$ and $g$
are increasing.
This concludes the proof of Theorem~\ref{thm:mono-Iso}.
\end{proof}    
  
\begin{prop}   \label{prop:1}
Let $a\geq 0$ and let $w_a:[0,1] \rightarrow \mathbb{R}$ be defined by
\[
w_a(x) 
=
\frac{\pFq{2}{1}{-a, -a}{1}{x}}{(x+1)^a}    
.
\]                                      
Then $w_a$ is: decreasing if $0 < a < 1$; increasing if $a > 1$; constant if $a \in\{ 0, 1\}$.
\end{prop}

\begin{proof}
Clearly, if $a \in\{ 0, 1\}$, then 
 $w_a(x)$ is constant, equal to~$1$ on $[0,1]$.

\smallskip Consider now the case $a>0$ with $a\neq 1$.
The derivative of $w_a(x)$ 
satisfies the hypergeometric identity
\begin{equation}\label{id:hyp}
\frac{w_a'(x) \cdot (x+1)^{a+1}}{a\cdot (a-1) \cdot (1-x)^{2a}}
\, = \, \pFq{2}{1}{a+1, a}{2}{x},
\end{equation}
which is a direct consequence of Euler's transformation formula~\cite[Eq. (2.2.7), p. 68]{AnAsRo99} and of Lemma~\ref{lem:1} with $a$ substituted by $-a$.

Since $a>0$, the right-hand side of~\eqref{id:hyp} has only positive Taylor coefficients, therefore it is positive on  $[0,1)$.
It follows that $w_a'(x) \geq 0$ on  $[0,1]$ if $a-1>0$,
and $w_a'(x) \leq 0$ on  $[0,1]$ if $a-1<0$.   
Equivalently, $w_a$ is increasing  on $[0,1]$ if $a>1$,
and decreasing on $[0,1]$ if $a<1$.
\end{proof}

\begin{lem}\label{lem:1}
  The following identity holds:
\[
(a+1)  (1-x) \cdot \, \pFq{2}{1}{a+1, a+2}{2}{x}
=
a (x+1) \cdot \, \pFq{2}{1}{a+1, a+1}{2}{x}    
+
\, \pFq{2}{1}{a, a}{1}{x}   .  
\] 
\end{lem}         

\begin{proof}
We will use two of the classical Gauss' contiguous relations~\cite[\S2.5]{AnAsRo99}:
\begin{equation}\label{cont:1}
\pFq{2}{1}{a+1, b+1}{c+1}{x}
=
\frac{c}{bx} \cdot \left( \pFq{2}{1}{a+1, b}{c}{x} - \pFq{2}{1}{a, b}{c}{x}     \right)
\end{equation}
and
\begin{align}\label{cont:2}
a\cdot & \left(
\pFq{2}{1}{a+1, b}{c}{x} -  \pFq{2}{1}{a, b}{c}{x}
\right) \nonumber
=\\
& \qquad \qquad \frac{(c-b) \cdot \pFq{2}{1}{a, b-1}{c}{x} + (b-c+ax) \cdot \pFq{2}{1}{a, b}{c}{x}     }{1-x} .
\end{align}
Applying~\eqref{cont:1} twice, once with $(b,c)=(a,1)$ and once with      $(b,c)=(a+1,1)$, 
the proof of the lemma is reduced to that of the identity
\[
(x-1) \cdot  \pFq{2}{1}{a+1, a+1}{1}{x}
+
2 \cdot  \pFq{2}{1}{a, a+1}{1}{x}
=
\pFq{2}{1}{a, a}{1}{x},
\]
which follows from~\eqref{cont:2} with $(b,c)=(a+1,1)$. 	                
\end{proof}

\smallskip \noindent 
\begin{rem}\label{rem:1}
A natural question is whether the function $\Iso$ enjoys higher monotonicity properties.
It can be easily seen that both $\Iso$ and its reciprocal $1/\Iso$ are neither convex nor concave. However we will prove that
 $z \mapsto \Iso(\sqrt{z})$ is concave and $z \mapsto 1/\Iso(\sqrt{z})$ is convex, on their  domain of definition $[0,3-2\sqrt{2})$. 

First recall that $1/\Iso(\sqrt{z}) = 
\frac{2^{5/4} \cdot \sqrt{\pi}}{3}
\cdot 
w_{1/2}(r(z))^{3/2} \cdot w_{3/2}(r(z))^{-1}$, 
where we set $r(z)= 4z/(1-z)^2$. 
Since $w_{1/2}=1/f$ and $w_{3/2}^{-1}=1/g$ are positive and decreasing, 
 while $r$ is nonnegative and increasing, 
proving that  $w_{1/2} \circ r$ and $w_{3/2}^{-1} \circ r$ are both convex
is enough to establish convexity of $z \mapsto 1/\Iso(\sqrt{z})$.

From $(\ref{id:hyp})$ and the chain rule it follows that 
\begin{equation}\label{id:war}
    \frac{\frac{\mathrm{d}}{\mathrm{d}z} \, w_a(r(z))}{4\, a \cdot (a-1) } = \pFq{2}{1}{a+1, a}{2}{\frac{4z}{(1-z)^2}} \cdot \frac{ (1-6\,z + z^2)^{2a}}{(1-z)^{4a}} \cdot \frac{(1-z)^{2a-1}}{(1+z)^{2a+1}}.
\end{equation}
We can justify convexity of both $w_{1/2}(r(z))$ and $w_{3/2}(r(z))^{-1}$ if we can prove that the right-hand side of (\ref{id:war}) is decreasing on $[0,3-2\sqrt{2})$. Moreover, it is easy to see that $(1-z)^{2a-1}/(1+z)^{2a+1}$ is decreasing on this interval for $a>3/2-\sqrt{2}$. Therefore, after changing variables $x={4z}/{(1-z)^2}$, it remains to show that 
\[
    \pFq{2}{1}{a+1, a}{2}{\frac{4z}{(1-z)^2}} \cdot \frac{ (1-6\,z + z^2)^{2a}}{(1-z)^{4a}} = \pFq{2}{1}{a+1, a}{2}{x} \cdot (1-x)^{2a}
\]
is decreasing for all $x \in [0,1)$. 
The derivative of the right-hand side is given by
\[
-\left( \frac{a(3-a)}{2} \cdot \pFq{2}{1}{a,a+1}{3}{x} + \frac{a(a+1)x}{6} \cdot \pFq{2}{1}{a+1,a+2}{4}{x} \right)\cdot (1-x)^{2a-1},
\]
hence is indeed negative for all $ x \in [0,1)$ if $0<a<3$. From this and
(\ref{id:war}) it follows that $1/\Iso(\sqrt{z}) $ is the product of two positive, decreasing and convex
functions and therefore inherits these properties. Finally, this also shows
that $\Iso(\sqrt{z})$ is both increasing and concave on $[0,3 - 2\sqrt{2})$.
\end{rem}

\medskip \noindent 
\begin{rem}\label{rem:2}
Bruno Salvy (private communication)
found an alternative short proof of Proposition~\ref{prop:1}.
The main idea is inspired by the Sturm–Liouville theory and the proof is based on the observation that $w_a(x)$ satisfies the linear differential equation (written in adjoint form):
\[
\frac{\mathrm{d}}{\mathrm{d}x} \left(x\left(\frac{1+x}{1-x}\right)^{2a} \cdot  \frac{\mathrm{d}}{\mathrm{d}x} w_a(x) \right)=  \frac{a(a-1)x }{ (1+x)^2}  \left(\frac{1+x}{1-x}\right)^{4a} \cdot w_a(x).
\]
The right-hand side is positive on $(0,1)$ when $a>1$ and negative if $0<a<1$. 
The same holds for its integral over $[0,t]$ for any $t<1$.
Looking at the left-hand side, 
this implies that $w_a'>0$ whenever $a>1$ and $w_a'<0$ when $0<a<1$.

We note that the same idea allows for a different proof of Remark~\ref{rem:1}. 
\end{rem}

\smallskip \noindent {\bf Acknowledgements.} 
We thank Mourad Ismail for his challenging question which led to Remark \ref{rem:1} and Bruno Salvy for his alternative proof of Proposition~\ref{prop:1}. Moreover, we are grateful to both for their interest in 
this work. The first author was supported in part by \textcolor{magenta}{\href{https://specfun.inria.fr/chyzak/DeRerumNatura/}{DeRerumNatura}} ANR-19-CE40-0018. The second author was supported by the \href{https://www.fwf.ac.at/}{Austrian Science Fund} (FWF): P-31338.

\end{document}